\documentclass{article}

\usepackage{amsfonts}
\usepackage{amsmath}
\usepackage{enumerate}

\newtheorem{theorem}{Theorem}

\newtheorem{corollary}[theorem]{Corollary}

\newtheorem{example}[theorem]{Example}

\newtheorem{lemma}[theorem]{Lemma}

\newtheorem{problem}[theorem]{Problem}

\newtheorem{remark}[theorem]{Remark}

\newcommand{\gl}{\lambda}

\newcommand{\gW}{\Omega}
\newcommand{\gD}{\Delta}
\newcommand{\gL}{\Lambda}
\newcommand{\mch}{\mathcal{H}}
\newcommand{\mcf}{\mathcal{F}}
\newcommand{\mbh}{\mathbb{H}}
\newcommand{\mbf}{\mathbb{F}}
\newcommand{\mbfl}{\mbf^L}
\newcommand{\mbr}{\mathbb{R}}
\newcommand{\mbd}{\mathbb{D}}
\newcommand{\mcd}{\mathcal{D}}
\newcommand{\mbn}{\mathbb{N}}
\newcommand{\mbg}{\mathbb{G}}
\newcommand{\mcg}{\mathcal{G}}

\newcommand{\mck}{\mathcal{K}}
\newcommand{\mbk}{\mathbb{K}}
\newenvironment{proof}[1][Proof]{\textbf{#1.} }{\ \rule{0.5em}{0.5em}}

\newenvironment{romenumerate}[1][0pt]{
\addtolength{\leftmargini}{#1}\begin{enumerate}
 }{\end{enumerate}}

\newcounter{thmenumerate}
\newenvironment{thmenumerate}
{\setcounter{thmenumerate}{0}%
 \def\item{\par
 \refstepcounter{thmenumerate}\textup{(\roman{thmenumerate})\enspace}}
}
{}

\newcommand{\refT}[1]{Theorem~\ref{#1}}
\newcommand{\refL}[1]{Lemma~\ref{#1}}
\newcommand{\refE}[1]{Example~\ref{#1}}
\newcommand{\tia}{totally inaccessible}

\newcommand\marginal[1]{\marginpar{\raggedright\parindent=0pt\tiny #1}}
\newcommand\SJ{\relax}
\newcommand\kolla{\relax}

\begingroup
  \divide\count255 by 60
  \multiply\count255 by -60
  \advance\count255 by \time
  \ifnum \count255 < 10 \xdef\klockan{\the\count1.0\the\count255}
  \else\xdef\klockan{\the\count1.\the\count255}\fi
\endgroup

\newcommand\set[1]{\ensuremath{\{#1\}}}

\newcommand\bigpar[1]{\bigl(#1\bigr)}

\def\rompar(#1){\textup(#1\textup)}    

\newcommand\E{\operatorname{\mathbb E{}}}

\newcommand\go{\omega}

\newcommand{\gU}{\Upsilon}

\newcommand\intoo{\int_0^\infty}
\newcommand\intot{\int_{0}^t}

\newcommand\REM[1]{{\raggedright\texttt{[#1]}\par\marginal{XXX}}}

\newcommand\ogl{{}^o\gl}

\newcommand\ett[1]{1_{\set{#1}}} 
\renewcommand\E{E}
\newcommand{\qs}{Q^\star}
\newcommand{\mcq}{\mathcal{Q}}

\begin{document}

\title{Absolutely Continuous Compensators}
\author{Svante Janson\thanks{
Uppsala University, 
Department of Mathematics, P.O. Box 480, 
SE-751 06 Uppsala, Sweden}\enspace  
and Sokhna M'Baye\thanks{
D\'epartement de de Math\'ematiques, 
\'Ecole Normale Sup\'erieure de Cachan, 
61 Avenue du Pr\'esident Wilson, 
  94235 Cachan Cedex, France}\enspace 
and Philip Protter\thanks{
School of Operations Research and Industrial Engineering, Cornell
University, Ithaca, New York 14853; Supported by NSF Grant  DMS-0906995} \\
\\
\\
}
\date{May 12, 2010} 
\maketitle

\begin{abstract}
\noindent We give sufficient conditions on the underlying filtration
such that all totally inaccessible stopping times have compensators
which are absolutely continuous.  If a semimartingale, strong Markov
process $X$ has a representation as a solution of a stochastic
differential equation driven by a Wiener process, Lebesgue measure,
and a Poisson random measure, then all compensators of totally
inaccessible stopping times are absolutely continuous with respect to
the minimal filtration generated by $X$.  However \c{C}inlar and Jacod
have shown that all semimartingale strong Markov processes, up to a
change of time and slightly of space, have such a representation. 
\end{abstract}

\section{Introduction}
\noindent  The celebrated Doob-Meyer Decomposition Theorem states that
if $X$ is a submartingale, then it can be written in the form $X=M+A$
where $X$ is a local martingale and $A$ is a unique, c\`adl\`ag 
increasing predictably measurable 
process with $A_0=0$.  (See, for
example,~\cite{PP}.) In the case of a point process of $N=(N_t)_{t\geq
  0}$ it is trivially a submartingale, and hence we know there exists
a process $A$ such that $N-A$ is a local martingale.  A special case
of interest in the theory of Credit Risk is the case 
\begin{equation}\label{ie1}
1_{\{t\geq R\}}-A_t=\text{ a martingale }
\end{equation}
The process $A$ in $(\ref{ie1})$ is known as \emph{the compensator of
  the stopping time $R$}, by an abuse of language.  It is common in
applications to assume, often without mention, that $A$ has absolutely
continuous paths.  That is, one often assumes \emph{a priori} that
$(\ref{ie1})$ is of the form 
\begin{equation}\label{ie2}
1_{\{t\geq R\}}-\int_0^t\gl_sds=\text{ a martingale }
\end{equation}
for some adapted process $\gl$.  The process $\gl$ is often referred
to as the \emph{hazard rate} and has intuitive content as the
instantaneous likelihood of the stopping time $R$ occurring in the
next infinitesimal time interval.  Of course this is not true in
general, and for example in the theory of credit risk K. Giesecke and
L. Goldberg have given a natural example where it does not 
hold~\cite[C3, p.~7]{GG}.  The goal of this paper is to give simple
and natural conditions on the generating underlying filtration to show
when the compensators of \emph{all} of the totally inaccessible
stopping times are absolutely continuous; that is, to give sufficient
conditions on the filtration such that they all have hazard rates.   

\section{Prior Results}
Previous work has been restricted to giving conditions on a given
stopping time in relation to the underlying filtration that ensures
the compensator is absolutely continuous.  Perhaps the most well known
of these conditions is that of S. Ethier and T.G. Kurtz~\cite{EK},
which we restate here.   
\begin{theorem}[Ethier--Kurtz Criterion]\label{t1} 
Let $\mbg=(\mcg_t)_{t\ge0}$ 
be a given filtration satisfying the usual hypotheses (see~\cite{PP} for the ``usual hypotheses.'')  Let $A$ be an increasing (not necessarily adapted) and integrable c\'adl\'ag 
process, with $A_0=0$. Let $\tilde{A}$ be the $\mbg$ compensator of A. If there is a constant $K$ such that for all $0\leq s\leq t$
\begin{equation}\label{ie3a}
E\{A_{t}-A_s\vert\mcg_s\}\leq K(t-s) \text{ a.s. }
\end{equation}
then the compensator of $A$ has absolutely continuous paths, a.s.  That is, it is of the form $\tilde{A}_t=\int_0^t\gl_sds$.
\end{theorem}
An extension of Theorem~\ref{t1} to necessary and
sufficient conditions for the $\mbg$ compensator to have such an
intensity process is given in the Cornell PhD thesis of Yan Zeng~\cite{Zeng}. A trivial extension is to replace the constant $K$ with an increasing predictable process $(K_t)_{t\geq 0}$, and then the inequality~(\ref{ie3a}) becomes:
\begin{equation}
E\{A_{t}-A_s\vert\mcg_s\}\leq K_s(t-s) \text{ a.s. }
\end{equation}
and of course the conclusion in Theorem~\ref{t1} still holds.
Zeng~\cite[p.~14]{Zeng} did a little better: 
\begin{theorem}[Yan Zeng]\label{t2} 
Let $A$ be an increasing (not necessarily adapted) and integrable measurable
process, with $A_0=0$. Let $\tilde{A}$ be the compensator of $A$. Then $d\tilde{A}_t\ll dt$ if and
only if there exists an increasing and integrable measurable process $D$ with $D_0=0$,
such that $d\tilde{D}_t\ll dt$  and for all $t\geq 0, h\geq 0$,
\begin{equation}
E\{A_{t+h}-A_t\vert\mcg_t\}\leq E\{D_{t+h}-D_t\vert\mcg_t\}
\end{equation}
and if equality holds, then we have $\tilde{A}=\tilde{D}$.
\end{theorem}

Another observation is perhaps useful to make.  Once a stopping time
has an absolutely  continuous 
compensator in a given filtration, say $\mbg$, then if it is also a
stopping time for a smaller filtration it also has an absolutely
continuous compensator in the smaller filtration.  Actually one can
obtain a more precise result, which is established in the book of
Martin Jacobsen~\cite{MJ}.  We provide here an original and elementary
proof of this result, and \emph{inter alia} we extend the result a little.   

\begin{theorem}\label{t3} 
Let $R$ be a $\mbg$ stopping time with compensator given by 
$\int_0^t\gl_sdc(s)$ for some  $\mbg$ adapted process $\gl$,
where $\mbg$ satisfies the usual hypotheses.  Here $s\mapsto c(s)$ is non random, continuous, and non-decreasing.  Let $\mbf$ be a
subfiltration of $\mbg$ also satisfying the usual hypotheses, and
suppose $R$ is also an $\mbf$ stopping time.  Then the $\mbf$
compensator of $R$ is given by $\int_0^t\,\ogl_sdc(s)$.\footnote{Note
  that for fixed $s$, $\ogl_s=E\{\gl_s\vert\mcf_s\}$ a.s.; the
  optional projection gives a method to define the projection via
  conditional expectation for all $s\geq 0$ simultaneously.  Since
  $\gl$ is positive, the optional projection exists.
  See~\cite{PP} for more details.}  That is we have 
\begin{align} 
&\text{If}& 1_{\{t\geq R\}}-\int_0^t\gl_sdc(s)\ &= \text{  a martingale in } \mbg \notag\\
&\text{then}& 1_{\{t\geq R\}}-\int_0^t\, \ogl_sdc(s)\ &= \text{  a martingale in } \mbf.
\end{align}
\end{theorem}
\begin{proof}
Let
\begin{equation}\label{e1a}
M_t=1_{\{t\geq R\}}-\int_0^t\gl_sdc(s).
\end{equation}
Then $M$ is a $\mbg$
martingale. 
Since $\gl_s\ge0$, the optional projection $\ogl_s$ exists with $0\le
\ogl_s\le\infty$. 
For every $s$, $\E\ogl_s=\E\gl_s$, and thus, by Fubini's theorem,
$\E\intoo\ogl_s\,dc(s)=\E\intoo\gl_s\,dc(s)=\E\ett{R<\infty}<\infty$. Thus, 
$A_t=\intot\ogl_s\,dc(s)$ is an integrable increasing continuous adapted
process; in particular, $\ogl_s<\infty$ for a.e.\ $s$ a.s.
We define 
\begin{equation}
  L_t=\ett{t\ge R} - A_t
=\ett{t\ge R} - \intot \ogl_s\,dc(s).
\end{equation}
If $0\le s\le t$ and $H$ is bounded and  $\mcf_s$ measurable, then, 
by Fubini's theorem and the fact that
for fixed $r$, $\ogl_r=\E\set{\gl_r\vert\mcf_r}$ a.s.,
\begin{equation}
  \begin{split}
\E\bigpar{H(L_t-L_s)}
&=\E\bigpar{H\ett{t\ge R>s}}
-\int_s^t \E \bigpar{H\E(\gl_r\vert\mcf_r)}\,dc(r)	
\\&
=\E\bigpar{H\ett{t\ge R>s}}
-\int_s^t \E (H \gl_r)\,dc(r)	
\\&
=\E\bigpar{H(M_t-M_s)}
=0.
  \end{split}
\end{equation}
Hence the uniformly integrable process $L_t$ 
is an $\mbf$ martingale.  This
gives 
the $\mbf$ canonical decomposition of the $\mbf$ submartingale
$1_{\{t\geq R\}}$ as 
$1_{\{t\geq R\}}=L_t+A_t$, and thus $A$ is the $\mbf$ compensator of
  $1_{\{t\geq R\}}$.  
\end{proof}

We include for emphasis the following obvious but important (and well known) corollary:
\begin{corollary} \label{nc1} 
Let $R$ be a $\mbg$ stopping time with compensator given by 
$\int_0^t\gl_sds$ for some  $\mbg$ adapted process $\gl$,
where $\mbg$ satisfies the usual hypotheses.  Let $\mbf$ be a
subfiltration of $\mbg$ also satisfying the usual hypotheses, and
suppose $R$ is also an $\mbf$ stopping time.  Then the $\mbf$
compensator of $R$ is given by $\int_0^t\,\ogl_sds$.  That is we have 
\begin{align} 
&\text{If}& 1_{\{t\geq R\}}-\int_0^t\gl_sds\ &= \text{  a martingale in } \mbg \notag\\
&\text{then}& 1_{\{t\geq R\}}-\int_0^t\, \ogl_sds\ &= \text{  a martingale in } \mbf.
\end{align}
\end{corollary}
Theorem~\ref{t3} and its Corollary~\ref{nc1} show that once there is a
filtration $\mbh$ such that a stopping time $R$ is totally inaccessible, if
$R$ has an AC compensator in $\mbh$, then it has an AC compensator in any
smaller filtration $\mbg$ as well.  In particular Dellacherie's result
(Theorem~\ref{th7} below) implies that \emph{the law of $R$ is absolutely
  continuous (ie, has a density) as well}.  We recall Dellacherie's result
here for the reader's convenience.  A proof can be found
in~\cite[p.\ 120]{PP}. 

\begin{theorem}[Dellacherie's Theorem]\label{th7}
Let $R$ be a nonnegative random variable with $P(R=0)=0, P(R>t)>0$ for each $t>0$.  Let $\mcf_t=\sigma(t\wedge R)$, the minimal filtration which renders $R$  a stopping time.  Let $F$ denote the law of $R$.  That is, $F(x)=P(R\leq x)$ for $x\geq 0$.  Then the compensator $A=(A_t)_{t\geq 0}$ of the process $1_{\{R\geq t\}}$ is given by 
$$
A_t=\int_0^t\frac{1}{1-F(u-)}dF(u).
$$
If $F$ is continuous, then $A$ is continuous, $R$ is totally inaccessible, and $A_t=-\ln (1-F(R\wedge t))$.
\end{theorem}

One may ask if that, once a compensator of a stopping time $R$ is a.s.\
singular with respect to Lebesgue measure, does that propagate down to
smaller filtrations, and in particular does it imply that the law of the
stopping time is singular as well?  The next example shows that this is not
true in general. 
\begin{example}\rm
Let $B$ be a standard one dimensional Brownian motion with natural filtration $\mbf$ and with a local time at zero $L=(L_t)_{t\geq 0}$.  Define the change of time
$$
\tau_t=\inf\{s>0: L_s>t\}.  
$$
Then $(\tau_t)_{t\geq 0}$ is a family of $\mbf$ stopping times.  Also,
$M_t=B_{\tau_t}$ is a local martingale for the filtration $\mbg$ given by
$\mcg_t=\mcf_{\tau_t}$ for $t\geq 0$.  Let $N$ be an independent Poisson
process, and consider the vector of processes on the appropriate product
space $(N_t-t,M_t)_{t\geq 0}$, with filtration $\mbh$ such that
$(N_t-t,M_t)$ is a vector of two martingales.  Then the family $(L_t)_{t\geq
  0}$ are stopping times for $\mbh$, and  
$$
N_{L_t}-L_t = \text{ a local martingale for the filtration } \tilde{\mbh}
$$
where $(\tilde{\mch}_t)=(\mch_{L_t})_{t\geq 0}$.  Since $L$ is Brownian local time at zero, it has paths which are singular with respect to Lebesgue measure, a.s.  However by Tanaka's formula we have
$$
\vert B_t\vert=\int_0^t\text{ sign}(B_s)dB_s+L_t
$$
and therefore $E(L_t)=E(\vert B_t\vert)=\sqrt{\frac{2}{\pi}}\sqrt{t}$.  

Next let 
$$
R=\inf\{s>0: N_{L_s}\geq 1\}.
$$
Then $1_{\{t\geq R\}}-L_{t\wedge R} =$ a martingale for the filtration
$\tilde{\mbh}$, and the compensator of $R$ 
is $L_{t\wedge R}$ which 
inherits the singular nature of
the paths of $L$.  That is, the compensator of $R$ has paths which are
a.s.\ singular with respect to Lebesgue measure.  However the law $F$ of $R$
satisfies 
$$
F(t)=P(R\leq t)=E(1_{\{R\leq t\}})=E(L_{t\wedge R}),
$$
which is absolutely continuous, since $t\mapsto
E(L_t)=\sqrt{\frac{2}{\pi}}\sqrt{t}$ is absolutely continuous.  Therefore by
Dellacherie's theorem the compensator of $R$ in the minimal filtration that
makes it a stopping time is absolutely continuous. 
\end{example}

\begin{corollary} Let $R$ be a stopping time for a filtration $\mbh$, and suppose there is a subfiltration $\mbg$ such that $R$ is totally inaccessible, and that the compensator of $1_{\{t\geq R\}}$ is given by $\int_0^t\gl_sdc(s)$ where $s\mapsto c(s)$ is non random, continuous, and non-decreasing. Let $F$ denote the law of $R$.  Then $dF(s)\ll dc(s)$.  In particular if $c(s)=s$, then $F$ is absolutely continuous; if $s\mapsto dc(s)$ is singular with respect to Lebesgue measure, then $F$ is also singular.
\end{corollary}

\begin{corollary} Let $R$ be a stopping time for a filtration $\mbh$.  Suppose there is a subfiltration $\mbg$ such that $R$ is totally inaccessible, and that the compensator of $1_{\{t\geq R\}}$ is given by $\int_0^t\gl_sdc(s)$ where $s\mapsto c(s)$ is non random, continuous, and non-decreasing. Then for any sub-subfiltration $\mbf$ where $R$ is still a stopping time, $R$ will still be totally inaccessible, and its compensator will be absolutely continuous with respect to $dc$.  In particular, a compensator cannot be singular of this form for a filtration $\mbg$ and then become absolutely continuous for a subfiltration $\mbf$.
\end{corollary}

\section{Filtration Level Results}
In their seminal paper of 1981, E. \c{C}inlar and \kolla
J. Jacod~\cite[Theorem 3.33 on page 206]{CJ} showed that any $\mbr^d$
valued strong Markov process which is a Hunt process, and which is
also a semimartingale, up to a change of time via an additive
functional ``clock,'' can be represented as the solution of a
stochastic differential equation driven by $dt, dW_t,$ and $n(ds,dz)$;
where $W$ is a standard multidimensional Brownian motion, and $n$ is a
standard Poisson random measure with mean measure given by
$ds\nu(dz)$.   

Therefore we assume as given a strong Markov Hunt process semimartingale which can be represented on a space $(\Omega,\mcf,\mbf,P^x)$ where $\mbf=(\mcf_t)_{t\geq 0}$, as follows:

\begin{eqnarray}\label{e1}
X_t=X_0&+&\int_0^tb(X_s)ds+\int_0^tc(X_s)dW_s\notag\\
&+&\int_0^t\int_{\mbr}k(X_{s-},z)1_{\{\vert k(X_{s-},z)\vert\leq 1\}}[n(ds,dz)-ds\nu(dz)]\notag\\
&+&\int_0^t\int_{\mbr}k(X_{s-},z)1_{\{\vert k(X_{s-},z)\vert>1\}}n(ds,dz)
\end{eqnarray}

\noindent We let $P^\mu$ denote the probability measure governing $X$ where the law of $X_0$ is $\mu$, and $\mbf^\mu$ denote the filtration containing $\mbf$ but such that $\mcf_0^\mu$ contains all of the $P^\mu$ null sets.  Our goal is to prove that the predictable compensator process $A$, of the process $1_{\{t\geq R\}}$ where $R$ is a given but arbitrary totally inaccessible stopping time, has absolutely continuous paths a.s.  The next  theorem makes this precise.

\begin{theorem}\label{t4}
For any totally inaccessible stopping time $R$ on the space $(\Omega,\mcf,\mbf^\mu,P^\mu)$ the predictable increasing process $A$, with $A_0=0$, such that $1_{\{t\geq R\}} -A_t=M_t$ is a martingale, has the form $A_t=\int_0^t\gl_sds$ for some adapted process $\gl$.
\end{theorem}

\begin{proof} 
It is now perhaps mostly forgotten, but most of the ingredients for the
proof of this theorem are contained in an old paper of P.A. Meyer, published 
in 1973~\cite{Meyer};
see also \cite{PAMIII}. 
The idea is to recall that for each law $P^\mu$, each
square integrable martingale $M$ with $M_0=0$ is null, if it is orthogonal
to all martingales of the form  
\begin{equation}\label{eq1}
C_t^n=g_n(X_t)-g_n(X_0)-\int_0^tLg_n(X_s)ds
\end{equation}
where $L$ is the infinitesimal  generator  of the underlying strong Markov
process $X$, and the functions $g_n$ are 
a suitable sequence of functions that are 
``nearly Borel'' measurable, belong
to the domain of $L$, and each $Lg_n$ is bounded.  Next let $D^n$ be the
sequence of (still square integrable) martingales obtained by an
orthogonalization procedure under $P^\mu$, for the sequence $C^n$.  Let $M$
denote an arbitrary and chosen square integrable martingale under $P^\mu$.
Then $M$ can be represented as a sum of stochastic integrals with respect
to the collection $D^n$, and hence $\langle M,M\rangle$ is absolutely
continuous with respect to the collection $\langle D^n,D^n\rangle$, hence
also absolutely continuous with respect to the collection $\langle
C^n,C^n\rangle$.  We next choose constants $\gl_n$ such that
$E(\sum_n\gl_n\langle C^n,C^n\rangle_t)<\infty$ for all $t>0$, and we let  
\begin{equation}
  \label{Kt}
K_t=\sum_n\gl_n\langle C^n,C^n\rangle_t.
\end{equation}
Thus, $\langle M,M\rangle$ is absolutely continuous with respect to $K$.

In the case of the martingale $M$ which is the compensated indicator
function of the stopping time $R$: 
$M_t=1_{\{R\geq t\}}-A_t$, with $A$ continuous
(which is equivalent to $R$ being totally inaccessible), then 
we have that 
$$
[M,M]_t=\sum_{s\leq t}(\Delta M_s)^2=1_{\{t\geq R\}}.
$$  
We conclude that  $\langle M,M\rangle_t=A_t$.
Thus $\langle M,M\rangle=A$ is also continuous, and we conclude that 
$\langle M,M\rangle$ is 
absolutely continuous with respect to the continuous part of $K$.  The
continuous part of $K$, however, is a version of the continuous additive
functional $H$ of the L\'evy system of $X$, as given in~(\ref{ec1}) which
follows this proof. 

Therefore $A=\langle M,M\rangle$ must be absolutely continuous with respect to $H$.  Finally for a Markov process
$X$ of the type given in~(\ref{e1}), we know that $dH_t$ is absolutely
continuous with respect to $dt$.  This completes the proof. 
\end{proof}

\noindent {\bf Remark:} Theorem~\ref{t4} gives a sufficient condition
for compensators of all totally inaccessible stopping times to be
absolutely continuous, within a semimartingale Hunt process framework.
One might ask for necessary and sufficient conditions.  The same proof
plus a use of L\'evy systems can provide this result, given in
Theorem~\ref{t5}.  The connection to L\'evy systems was recently
recalled in the work of X. Guo and Y. Zeng~\cite{GZ}, and examples of
intensities arising in the field of Credit Risk can be found there and
in their references, as well as in~\cite{GG} and~\cite{JPS}, for
example.  See also~\cite{JP}. Examples of intensities arising in the
field of Survival Analysis can be found in the book of Fleming and
Harrington~\cite{FH}. 

\begin{corollary}\label{t5} 
Let $X=(\gW,X,P^\mu)$ be a semimartingale Hunt process with a L\'evy system $(K,H)$, where $K$is  a kernel on $\mbr$ and $H$is  a continuous additive functional, given by the following relationship:
\begin{equation}\label{ec1}
E^\mu\left (\sum_{0<s\leq t}f(X_{s-},X_{s})1_{\{X_{s-}\neq X_s\}}\right)=E^\mu\left(\int_0^tdH_s\int_\mbr K(X_s,dy)f(X_s,y)\right)
\end{equation}
Then all totally inaccessible stopping times have absolutely continuous
compensators if and only if the continuous additive functional of
equation~\eqref{ec1} is absolutely continuous with respect to Lebesgue
measure, a.s.  That is, if and only if $dH_s\ll ds$ a.s. 
\end{corollary}
\begin{proof}
The representation of a semimartingale Markov process given in
equation~(\ref{e1}) assumes there has already been a time change, if
necessary, to arrive at a Poisson random measure with compensator
$ds\,\nu(dx)$.  Here we are not making that assumption.  The results contained
in (for example)~\cite{Meyer} 
and \cite{BenvJac} 
show that for the additive functional $H$ of
the L\'evy system,  any representation such as~(\ref{e1}) must have that the
compensator of the corresponding ``Poisson random measure" will be
absolutely continuous in the $t$ variable with respect to $H$.  The additive
functional $H$ is not necessarily unique within the framework of L\'evy
systems, but any other version will be mutually absolutely continuous with
respect to it. 
Therefore by the proof of Theorem~\ref{t4} we have that all
totally inaccessible times are absolutely continuous with respect to $dt$ if
$dH_t\ll dt$.

For the necessity, suppose 
that every totally inaccessible stopping time has absolutely continuous
compensator.
Since $C^n$ in \eqref{eq1} jumps only when $X$ jumps, and the jumps
of the Hunt process $X$ can be covered by a countable collection of totally
inaccessible stopping times, it follows that $d\langle C^n,C^n\rangle_t\ll
dt$. Hence, by \eqref{Kt}, $d K_t\ll dt$. In particular $K_t$ is continuous
and is thus a version of the additive functional $H$, so $dH_t\ll
dK_t\ll dt$. 
\end{proof}

A  useful result related to Theorem~\ref{t4} is the following.  Jacod and Skorohod~\cite{JacodS} define a
\emph{jumping filtration} $\mbf$ to be a filtration such that there
exists a sequence of stopping times $(T_n)_{n=0,1,\dots}$ increasing
to $\infty$ a.s.\ with $T_0=0$ and such that for all $n\in\mbn, t>0$,
the $\sigma$-fields $\mcf_t$ and $\mcf_{T_n}$ coincide on $\{T_n\leq
t<T_{n+1}\}$.  We then have: 

\begin{theorem}\label{t5a}
Let $N=(N_t)_{t\geq 0}$ be a point process without explosions that generates a quasi-left continuous jumping filtration, and suppose there exists a process $(\gl_s)_{s\geq 0}$ such that 
\begin{equation}
N_t-\int_0^t\gl_sds= \text{  a martingale. }
\end{equation}
Let $\mbd=(\mcd_t)_{t\geq 0}$ be the (automatically right continuous)
filtration generated by $N$ and completed in the usual way.  
Then for any $\mbd$ totally inaccessible stopping time $R$ we have that the
compensator of $1_{\{t\geq R\}}$ has absolutely continuous paths,
a.s. 
\end{theorem}
\begin{proof}
\noindent By Theorem 2 of~\cite{JacodS} we know that
$\{R<\infty\}\subset \bigcup_{n\geq 1}\{R=T_n\}$, a.s.  
This implies that
$1_{\{t\geq R\}}\leq N_t$.  We write 
\begin{equation}\label{e2}
N_t=1_{\{t\geq R\}}+(N_t-1_{\{t\geq R\}})=1_{\{t\geq R\}}+C_t
\end{equation}
Let us now take compensators of both sides of (\ref{e2}), and let $A$ denote the compensator of $1_{\{t\geq R\}}$, and $\tilde{C}$ denote the compensator of $C$.  We want to show $dA_t\ll dt$.  Then the compensators version of equation~(\ref{e2}) becomes 
\begin{equation}
\int_0^t\gl_sds=A_t+\tilde{C}_t
\end{equation}
since taking compensators is a linear operation.  Since both $dA_t$ and $d\tilde{C}_t$ are positive measures on $\mbr_+$, it follows that $dA_t\ll \gl_tdt$ and $d\tilde{C}_t\ll \gl_tdt$.
\end{proof}

\begin{corollary}\label{c1}
Let $N$ be a Poisson process with parameter $\gl$, and $R$ be a
totally inaccessible stopping time on the minimal space generated by
$N$.  Then the compensator\footnote{The compensator of a stopping time
$R$ refers to the compensator of the increasing process $1_{\{t\geq
R\}}$.} of $R$ has paths which are absolutely continuous. 
\end{corollary}

A result which is related to Theorem~\ref{t5a}, but does not involve a
hypothesis on the filtration, is the following. 
For convenience we define $\gD Z_\infty=0$; 
hence $\set{\gD Z_R>0}\subseteq\set{R<\infty}$.

\begin{theorem}\label{t6}
Suppose $Z$ is an increasing process which has an absolutely
continuous compensator; that is, suppose there exists a nonnegative
adapted process $\gl$ such that $Z_t-\int_0^t\gl_sds=$ a martingale.
Let $R$ be a stopping time such that 
$P(\gD Z_R>0\cap\{R<\infty\})=P(R<\infty)$,
i.e., $\gD Z_R>0$ a.s.\ on \set{R<\infty}. 
Then $R$ has an absolutely
continuous compensator.  That is, there exists a nonnegative adapted
process $\mu$ such that $1_{\{t\geq R\}}-\int_0^t\mu_sds=$ a
martingale.   
\end{theorem}

\begin{proof} 
Let  
\begin{eqnarray}
\gL_1&=&\{\gD Z_R\geq 1\}, \notag\\
\gL_n&=&\{\frac{1}{n}\leq\gD Z_R<\frac{1}{n-1}\}, \qquad n\ge2,\notag\\
Y^n_t&=&1_{\{t\geq R_{\gL_n}\}.}
\end{eqnarray}
We have that $nZ_t-Y^n_t$ is again an increasing process, and if we
observe that 
\begin{equation}\label{e3a}
nZ_t=(nZ_t-Y^n_t)+Y^n_t,
\end{equation}
then by taking compensators of both sides in~(\ref{e3a}) we have, by
the same argument as in the proof of Theorem~\ref{t5}, that the
compensator of $Y^n$ is absolutely continuous.  Therefore we can write 
the compensator as $\int_0^t\mu^n_sds$, so 
\begin{equation}\label{e4a}
Y^n_t-\int_0^t\mu^n_sds = \text{ a martingale}.
\end{equation}
Furthermore, $1_{\{t\geq R\}}=\sum_{n=1}^\infty Y_t^n$   
and thus the compensator $A$ of $1_{\{t\geq R\}}$ is
\begin{equation}
A_t=\sum_{n=1}^\infty\int_0^t\mu^n_sds
=\int_0^t\left(\sum_{n=1}^\infty\mu_s^n\right)ds
\end{equation}
by the Fubini--Tonelli theorem, and the theorem is proved.
\end{proof}

\section{Related Results}

In this section we relate the results of the preceding part of the paper to some situations that arise in Mathematical Finance Theory.  Indeed, it is often the case in Finance that one begins with a system $(\Omega,\mcf,P,\mbf,X)$ and then one changes to an equivalent probability measure $Q$ such that $X$ is a sigma martingale (or less generally, a local martingale) under $Q$.  We will say that a probability measure $Q$ has {\bf Property AC} if all totally inaccessible stopping times have absolutely continuous compensators under $Q$.

\begin{theorem}\label{nt1}
 Suppose that $(\Omega,\mcf,P,\mbf,X)$ is a given system, and that
 there exists an equivalent probability measure $Q^\star$ (which can
 be $P$ itself) such that $\qs$ has Property AC. If $\mcq$ is the
 set of all  probability measures equivalent to $P$, we have that
 Property AC holds under any $Q\in\mcq$.  
\end{theorem}

\noindent\begin{proof}
Suppose $\qs\in\mcq$, and let $R\in\mcq$, so that $R$ is equivalent to
$\qs$. Let $\tau$ be a totally inaccessible stopping time, so that we
can write  
$$
1_{\{t\geq\tau\}}-\int_0^t\gl_sds=\text{ a martingale, under } \qs.
$$
By the predictable version of the Meyer--Girsanov theorem (see, eg,~\cite[p.\ 135]{PP})
\begin{eqnarray}\label{ne1}
Z&=&\frac{dR}{d\qs}\notag\\
Z_t&=&E_{\qs}\{\frac{dR}{d\qs}\vert\mcf_t\}\notag\\
M_t&=&1_{\{t\geq\tau\}}-\int_0^t\gl_sds\notag\\
M_t&=&\left(1_{\{t\geq\tau\}}-\int_0^t\gl_sds-\int_0^t\frac{1}{Z_{s-}}
d\langle Z,M\rangle_s\right)+\int_0^t\frac{1}{Z_{s-}}d\langle Z,M\rangle_s
\end{eqnarray}
and we have the term in parentheses on the right side of~(\ref{ne1})
is a martingale under $R$.  
Therefore the compensator of $\tau$ under $R$ is 
$\int_0^t\gl_sds+\int_0^t\frac{1}{Z_{s-}}d\langle Z,M\rangle_s
$,
and is absolutely
continuous as soon as $d\langle Z,M\rangle_t\ll dt$.  We next note
that  
$$
[M,M]_t=\sum_{s\leq t}(\Delta M_s)^2=1_{\{t\geq\tau\}}
$$
since $M$ has only one jump, and it is of size one.  Since
compensators are unique, this means that $\langle
M,M\rangle_t=\int_0^t\gl_sds$, since $\langle M,M\rangle_t$ is the
compensator of $[M,M]$, and it exists since $[M,M]_t\in L^1$.
Moreover we also know that $\langle Z,M\rangle_t$ exists, and by the
Kunita--Watanabe inequality (see, eg,~\cite[p.\ 150]{PP}) we have that
$d\langle Z,M\rangle_t\ll d\langle M,M\rangle_t\ll dt$, a.s.  The
result follows. 
\end{proof} 

A topic that has achieved importance recently in the theory of Credit
Risk is that of the Expansion of Filtrations.  See, for
example,~\cite{B-SJ},~\cite{CN}, or~\cite{JLC}.  In the case of initial
expansions, we can expand using Jacod's criterion (see~\cite{Jacod}
or~\cite[p.\ 371]{PP}) by adding a random variable $L$ to the
filtration $\mbf$ at time $0$, provided that for each $t\geq 0$ the
(regular) conditional distribution of $L$ given $\mcf_t$, denoted
$\gU_t(\go,dx)$, is such that $\gU_t(\go,dx)\ll \eta_t(dx)$, where
$\eta_t(dx)$ is a  $\sigma$-finite measure.  (The key part is that $\eta_t(dx)$ does not depend on
$\go$.) It is shown~(\cite{Jacod},~\cite{PP}) that we can in general
replace $\eta_t(dx)$ with $\eta(dx)$ which does not depend on $t$.  We
define $q^x_t$ to be a c\`adl\`ag martingale such that
$\gU_t(\go,dx)=q^x_t\eta(dx)$. Finally, we let
$\mcg_t^0=\mcf_t\vee\sigma(t\wedge L)$, and $(\mcg_t)_{t\geq 0}$ be
the completed, right continuous version of the filtration $\mbg^0$.  

We return to considering the probability measure $P$ only; the results 
generalize immediately to all measures $\qs\in\mcq$ by Theorem \ref{nt1}.

\begin{theorem}\label{nt2}
Suppose we expand the filtration $\mbf$ by adding a random variable
$L$ initially, where its conditional distribution $\gU_t(\go,dx)\ll
\eta_t(dx)$ for some $\sigma$-finite measure $\eta_t(dx)$. \SJ
 Suppose also that $P$ has
Property AC.  Then $P$ has Property AC under the expanded filtration
$\mbg$.
\end{theorem}

\noindent\begin{proof}
Let $\tau$ be a totally inaccessible stopping time and 
recall that $M_t= 1_{\{t\geq\tau\}}-\int_0^t\gl_sds$ is our decomposition
 for the filtration $\mbf$.  
Jacod's theorem \cite[Th\'eor\`eme (2.5)]{Jacod}  \kolla
gives that, for some predictable process $k^x_s$  
and with $\langle\cdot,\cdot\rangle$ computed always in $\mbf$,
\begin{eqnarray}\label{ne2}
\langle q^x,M\rangle_t&=&\int_0^tk^x_sq^x_{s-}d\langle M,M\rangle_s,
\quad\text{when $q^x_{t-}>0$, for $\eta$-a.e.\ $x$},
 \notag\\
M_t-\int_0^tk_s^Ld\langle M,M\rangle_s&=&\text{ a martingale in } \mbg.
\end{eqnarray}
Since $M_t= 1_{\{t\geq\tau\}}-\int_0^t\gl_sds$  in $\mbf$, we have that the
compensator of $1_{\{t\geq\tau\}}$ in $\mbg$  is again absolutely
continuous, since as we saw in the proof of Theorem~\ref{nt1}, $d\langle
M,M\rangle_t$  is absolutely continuous.  
\end{proof}

We have an analogous result for progressive expansions.  Using the notation
and results presented in~\cite[Chapter VI, Section 3]{PP}, we let $L$ be a
positive random variable that avoids all $\mbf$ stopping times.  That is,
$P(L=\nu)=0$ for all $\mbf$ stopping times $\nu$.  Since constants are
stopping times, we note that this implies that $L$ has a continuous
distribution function.  By Dellacherie's theorem (Theorem~\ref{th7}) we have
that this implies $L$ is totally inaccessible, at least in the minimal
filtration that turns $L$ into a stopping time;
see further \refL{LL} below.

We let
$Z_t={}^o1_{\{L>t\}}$, 
where the $o$ superscript denotes optional projection
onto the $\mbf$ filtration.  We let $A^L=(A^L_t)_{t\geq 0}$ denote the
(predictable) compensator of $1_{\{t\geq L\}}$ for the filtration $\mbf$.
(The process $1_{\{t\geq L\}}$ need not be adapted in order to have a
compensator.)   We then have that the Doob--Meyer decomposition of $Z$ in
$\mbf$ is $Z_t=M^L_t-A^L_t$, where $M^L_t$ is defined by $M^L_t\equiv
Z_t+A^L_t$.  

Recall 
that a random time $L$ is called an {\bf honest time} if it is the end of an
optional set.  More precisely a random variable $L$ is called {\bf honest}
if for every $t\leq \infty$ there exists an $\mcf_t$ measurable random
variable $L_t$ such that $L=L_t$ on $\{L\leq t\}$.  (See, e.g.,~\cite[p.\
  381--382]{PP} for more on honest times.) 
If $L$ is honest, we let $\mbf^L$ be the filtration defined by  
\begin{equation}\label{fl}
\mcf^L_t=\{\Gamma:\Gamma=(A\cap\{L>t\})\cup(B\cap\{L\leq t\})\text{ for some }A,B\in\mcf_t\}.
\end{equation}
(It is easy to see that $\mbf^L$ is a filtration satisfying the usual  
hypotheses when $L$ is honest, see \cite[Theorem VI.17]{PP}.)
We note that this definition of $\mbf^L$ is not the standard one.
In~\cite{PP} it is called the filtration $\mbg$. The standard definition of
$\mbf^L$ (which does not require $L$ to be honest) 
is as follows, and to avoid confusion, we denote it $\mbk^L$:  
$$
\mck^L_t=\{\Gamma\in\mcf:\text{there exists }\Gamma_t\in\mcf_t: \Gamma\cap\{L>t\}=\Gamma_t\cap\{L>t\}\}
$$
(Thus $\mcf^L_t\subseteq\mck^L_t$.) 
The filtration $\mbf^L$ is called the {\bf progressive expansion} of
$\mbf$ under $L$.
We assume that $L$ is honest, and consider the filtration $\mbf^L$. It is
easy to see that  $L$ becomes a stopping time for $\mbf^L$.  
($\mbf^L$ is
the smallest expansion of $\mbf$ that makes $L$ a stopping time.) 

\begin{theorem}\label{th8}
Let $L$ be a positive honest random variable which avoids all $\mbf$
stopping times.  If $P$ has Property AC for $\mbf$, then for any totally
inaccessible $\mbf$ stopping time $\tau$, the  compensator of
$\tau$ in $\mbf^L$ is again absolutely continuous. 
\end{theorem}  

\begin{proof}
We begin by assuming that $\tau$ is a totally inaccessible $\mbf$ stopping
time.  It is shown in~\cite[Theorem~VI.18]{PP} 
that for a square integrable martingale
$X$, 
its decomposition  under $\mbf^L$, is given by  
\begin{eqnarray*}
X_{t}&=&
\left(X_{t}-\int_0^{t\wedge L}\frac{1}{Z_{s-}}d\langle X,M^L\rangle_s+1_{\{t\geq L\}}\int_L^t\frac{1}{1-Z_{s-}}d\langle X,M^L\rangle_s\right)\\
&&+\left(\int_0^{t\wedge L}\frac{1}{Z_{s-}}d\langle X,M^L\rangle_s-1_{\{t\geq L\}}\int_L^t\frac{1}{1-Z_{s-}}d\langle X,M^L\rangle_s\right).
\end{eqnarray*}
In our case, the $\mbf$ martingale $X$ is equal to $M$, where
$M_t=1_{\{t\geq\tau\}}-\int_0^t\gl_sds$.  But we already know that $d\langle
M,M\rangle_t\ll dt$, whence again by the Kunita--Watanabe inequality we have
that $d\langle X,M^L\rangle_t\ll d\langle X,X\rangle_t= d\langle
M,M\rangle_t\ll dt$, and the result follows for all totally inaccessible
stopping times $\tau$.  
\end{proof}

\refT{th8} shows the AC property only for $\mbf$ stopping times
$\tau$. In order to extend it to $\mbf^L$ stopping times, we need two
lemmas (and an extra condition). Note that the first part of \refL{L1}(i) is a special case of the
Lemma in \cite[p.~378]{PP} (with essentially the same proof).

\begin{lemma}\label{L1}
Let $T$ be a stopping time for $\mbf^L$.  Then:
\begin{romenumerate}
  \item
There exists an $\mbf$ stopping time $S$ such that $T\wedge L=S\wedge L$
a.s. 
If\/ $T$ is $\mbf^L$ \tia, then $S$ can be chosen to be $\mbf$ \tia.
  \item
There exists a sequence $(S_n)_{n\ge1}$ of\/ $\mbf$ stopping times such that
$[T]\subseteq[L]\cup\bigcup_{n=1}^\infty[S_n]$ a.s.
If\/ $T$ is $\mbf^L$ \tia, then all $S_n$ can be chosen to be $\mbf$ \tia.
\end{romenumerate}
\end{lemma}

\begin{proof}
The process $U_t=\ett{t>T}$ is $\mbf^L$ predictable, so by \cite[Theorem
  VI.17]{PP}, there exist two $\mbf$ predictable processes $H$ and $K$ such
  that 
$$
U=H1_{[0,L]}+K1_{(L,\infty)}.
$$
Define 
\begin{align*}
  R_0=\inf\set{t\ge0:H_t\neq0}
\intertext{and, for rational $r>0$,}
  R_r=\inf\set{t>r:K_t\neq0}.
\end{align*}
These are $\mbf$ stopping times. 
If $T<L$, then $H=0$ on $[0,T]$ and $H=1$ on $(T,L]$, so $R_0=T$.
If $T>L$, then $K=0$ on $(L,T]$ and $K=1$ on $(T,\infty)$, so $R_r=T$ for
  every $r\in(L,T)$.
Consequently,
$[T]\subseteq[L]\cup\bigcup_{r\ge0}[R_r]$,
and the first assertion in (ii) follows by rearranging the countable set
of stopping times $(R_r)_{r\ge0}$ into a sequence $(S_n)_{n\ge1}$.

Furthermore, if $T\ge L$, then $R_0\ge L$, and thus $T\wedge L=R_0\wedge L$,
so the first part of (i) follows with $S=R_0$. 

Now suppose that $T$ is \tia. Each $\mbf$ stopping time $R_r$ can be
decomposed into two $\mbf$ stopping times as $[R_r]=[R^a_r]\cup[R^i_r]$ with $R^a_r$ accessible and $R^i_r$ \tia\  \cite [p.~104]{PP}.
Then   $R^a_r$ is also for $\mbf^L$ an accessible stopping time, and since
$T$ is \tia, $P(T=R^a_r)=0$. Hence
$[T]\subseteq[L]\cup\bigcup_{r\ge0}[R^i_r]$ a.s.\
and we can replace $R_r$ by $R^i_r$ above.
\end{proof}

\begin{lemma}
  \label{LL}
$L$ is a totally inaccessible stopping time for $\mbf^L$.
\end{lemma}
\begin{proof}
$L$ is a stopping time by the definition of $\mbf^L$.

Suppose that $T$ is a $\mbf^L$ predictable stopping time, and let $T_n$ be a 
sequence of $\mbf^L$ stopping times that announces $T$, i.e., $T_n\nearrow T$ and
$T_n<T$ for all $n$ when $T>0$ \cite[Section III.2]{PP}.
By \refL{L1}, there exist $\mbf$ stopping times $S_n$ such that $T_n\wedge
L=S_n\wedge L$. 
Let $S=\liminf S_n$; this is an $\mbf$ stopping time. 
On the set $\{0<T\le L\}$, each $T_n<T\le L$, so $T_n\wedge L=T_n$ and
$S_n=T_n$; thus $S=T$. Further, on $\{T=0\}$, trivially each $T_n=0<L$ a.s.,
so $S_n=0$ and $S=0=T$ a.s. 
Hence, $S=T$ a.s.\ on $\{T\le L\}$.
Consequently,
$$
P(T=L)=P(S=T=L)\le P(S=L)  =0,
$$
because $S$ is an $\mbf$ stopping time.
Since $T$ is an arbitrary $\mbf^L$ predictable stopping time, this shows
that $L$ is \tia.
\end{proof}

We can now extend Theorem~\ref{th8} to $\mbf^L$ totally inaccessible
stopping times, but we need an extra condition.

\begin{theorem}\label{th9} 
Let $L$ be a positive honest random variable which avoids all\/ $\mbf$
stopping times, and suppose that $P$ has Property AC for $\mbf$.
Then $P$ has Property AC for $\mbf^L$
if and only if the compensator $A^L$ of $L$ in 
$\mbf$ is absolutely continuous on $[0,L]$. 
\end{theorem}

\begin{proof}
  By the Jeulin--Yor theorem \cite{JeulinYor,GZ},
the $\mbf^L$ compensator of $L$ is 
$\int_0^{t\wedge L}\frac1{Z_{s-}} dA^L_s$. 
Hence this compensator is absolutely continuous if and only if $A^L_t$ is
absolutely continuous on $[0,L]$.

The necessity of this condition is now clear, since $L$ is \tia\ by
\refL{LL}.

Conversely, suppose that this condition holds so that $L$ has an absolutely
continuous compensator. Let $T$ be a \tia\ stopping time for $\mbf^L$. 
By \refL{L1}, there exists a sequence $(S_n)_{n\ge1}$ of \tia\ $\mbf$
stopping times such that $[T]\subseteq[L]\cup\bigcup_{n=1}^\infty[S_n]$.
For notational convenience, let $S_0=L$, so
$[T]\subseteq\bigcup_{n=0}^\infty[S_n]$.

By assumption (for $n=0$) and \refT{th8} (for $n\ge1$), each $S_n$ has an
absolutely continuous compensator $A^n_L$ for $\mbf^L$; we write
$A^n_t=\int_0^t\lambda^n_s\, ds$. Let $T_n=T$ if $T=S_k$ for some $k\le n$,
and $T_n=\infty$ otherwise. Then $T_n$ is a stopping time with
$[T_n]\subseteq\bigcup_{k=0}^n[S_k]$, and it follows that $T_n$ has
a compensator $B^n$ for $\mbf^L$ 
such that $\sum_{k=0}^n A^k_t-B^n_t$ is an increasing
  process; thus $B^n_t=\int_0^t f^n_s\,ds$ with $0\le
  f^n_s\le\sum_{k=0}^n\lambda^k_s$; in particular the compensator $B^n$ of
  $T^n$ is  absolutely continuous. It now follows by monotone convergence
  that $B^n$ converges to the compensator $B$ of $T$ for $\mbf^L$, and thus
  this  compensator too is absolutely continuous. 
\end{proof}

The condition that $A^L$ be absolutely continuous on $[0,L]$ \SJ
seems, unfortunately, to be rather restrictive. As we see in \refE{Eexit}
below, in several natural examples, $A^L$ is, on the contrary, singular,
which by the proof above implies that the $\mbfl$ compensator of $L$ is
singular, and thus that $P$ does \emph{not} have Property AC for $\mbfl$.
Nevertheless, \refT{th8} still applies, and more generally, the $\mbfl$
compensator of every \tia\ $\mbfl$ stopping time $T$ such that $P(T=L)=0$ is
absolutely continuous.

\begin{problem}\label{p1}
  \begin{thmenumerate}
\item
Does there exist an honest time $L$ (for some $(\Omega,\mcf,P,\mbf)$)
such that $A^L$ is absolutely continuous a.s.? 	
\item
In particular, does there exist such an honest time for the natural
filtration of a standard Brownian motion?
  \end{thmenumerate}
\end{problem}
We note that the relatively recent work of A. Nikeghbali~\cite{AN} makes the positive resolution of Problem~\ref{p1} seem unlikely.
\begin{remark}\rm
It is easily seen that $A^L$ always is continuous, since otherwise the set
of jump times would be an $\mbf$ predictable set and thus there would exist
 a   predictable $\mbf$ stopping time $R$ such that $P(\Delta
 A^L_R>0)>0$. But 
  then 
 $E (\Delta \ett{t\ge L})_R=E\Delta A^L_R>0$ so $P(L=R)>0$, a contradiction.
\end{remark}

\begin{example}\rm  \label{Eexit}  \SJ
Typical examples of honest times are exit times. For a simple example,
consider a standard Brownian motion $B_t$ with its
  standard filtration $\mbf$, and let $L=\sup\{t\le1:B_t=0\}$. Then $L$ is
  an honest time, and $P(L=T)=0$ for every stopping time $T$ (by the strong
  Markov property of $B$, which implies that a.s.\ 
either $B_T\neq0$ \SJ
or $B_{T+t}=0$ for some
  sequence of $t\searrow0$).
Since $L$ belongs to the predictable set $\{t:B_t=0\}$, the compensator
  $dA^L$ is a.s.\ supported by this set, but this set has Lebesgue measure
  $0$, so $dA^L$ is a.s.\ singular.

In fact \cite{AzJKY,YorII}, a simple calculation shows that
for $t<1$, 
$Z_t=P(L>t|\mcf_t)=2\Phi(-|B_t|/\sqrt{1-t})$, where $\Phi$ is the standard
  normal distribution function, and 
$dA^L_t = \sqrt{\frac{2}{\pi(1-t)}} dL^0_t$,
where $L^0$ is the local time at $0$.

Several similar examples with singular $A^L$ are treated in \cite{AzJKY} and
\cite{YorII},  
for example
$\sup\set{t<1:B_t\in F}$ for a given finite set $F$,
$\sup\set{t<\tau_1:B_t=0}$ with
$\tau_1=\inf\set{t>0:B_t=1}$,
and
$\sup\set{t:|B^{(3)}_t|\ge1}$ where
$B^{(3)}$ is a three-dimensional Brownian motion
and thus $|B^{(3)}|$ is a $\mathrm{BES}(3)$ process
 (see also \cite{Getoor} for this exit time).
\end{example}

The results above are easily extended by induction to the case of 
the filtration $\mbf^{\{L^1,L^2,\dots\}}$ extended
by a finite or infinite, strictly increasing sequence  
$(L_n)_{n=1}^N$ of positive honest variables that avoid all $\mbf$
stopping times; cf.\ \cite[p.\ 190--191]{DM}. We omit the details.
One could also consider a more general setting, where for
example $L^{i+1}$ is honest only for $\mbf^{L^i}$, 
instead of requiring it
and all subsequent random times to be honest for $\mbf$, 
or the case where the $L^i$ need not be required to be strictly
increasing.  We do not treat these cases here.  The necessary theory to do
so is available, however, within the book of Th.\ Jeulin~\cite{Jeulin}. 

\end{document}